\documentclass[11pt]{article}

\usepackage{latexsym}
\usepackage{amssymb}

\newtheorem{thm}{Theorem}

\newtheorem{lem}{Lemma}

\begin{document}
{
\begin{center}
{\Large\bf
On extensions of $J$-skew-symmetric and $J$-isometric operators.}
\end{center}
\begin{center}
{\bf S.M. Zagorodnyuk}
\end{center}

\section{Introduction.}
Last years an increasing number of papers was devoted to the investigations of operators related to
a conjugation in a Hilbert space, see, e.g. ~\cite{cit_100_GP}, \cite{cit_200_GP2}, \cite{cit_300_Z}, \cite{cit_400_LZ} and references therein.
A conjugation $J$ in a Hilbert space $H$ is an {\it antilinear} operator on $H$ such that $J^2 x = x$, $x\in H$,
and $(Jx,Jy)_H = (y,x)_H$, $x,y\in H$.
The conjugation $J$ generates the following bilinear form:
$$ [x,y]_J := (x,Jy)_H,\qquad x,y\in H. $$
A linear operator $A$ in $H$ is said to be $J$-symmetric ($J$-skew-symmetric) if
\begin{equation}
\label{f1_2}
[Ax,y]_J = [x,Ay]_J,\qquad x,y\in D(A),
\end{equation}
or, respectively,
\begin{equation}
\label{f1_3}
[Ax,y]_J = -[x,Ay]_J,\qquad x,y\in D(A).
\end{equation}
A linear operator $A$ in $H$ is said to be $J$-isometric if
\begin{equation}
\label{f1_3_1}
[Ax,Ay]_J = [x,y]_J,\qquad x,y\in D(A).
\end{equation}
If $\overline{D(A)} = H$, then conditions~(\ref{f1_2}), (\ref{f1_3}) and (\ref{f1_3_1}) are equivalent to the following conditions:
\begin{equation}
\label{f1_3_2}
JAJ \subseteq A^*,
\end{equation}
\begin{equation}
\label{f1_3_3}
JAJ \subseteq -A^*,
\end{equation}
and
\begin{equation}
\label{f1_3_4}
JA^{-1}J \subseteq A^*,
\end{equation}
respectively.
A linear operator $A$ in $H$ is called $J$-self-adjoint ($J$-skew-self-adjoint, or $J$-unitary) if
\begin{equation}
\label{f1_4}
JAJ = A^*,
\end{equation}
\begin{equation}
\label{f1_5}
JAJ = -A^*,
\end{equation}
or
\begin{equation}
\label{f1_6}
JA^{-1}J = A^*,
\end{equation}
respectively.

We shall prove that each densely defined $J$-skew-symmetric operator
(each $J$-isometric operator with $\overline{D(A)}=\overline{R(A)}=H$) in a Hilbert space $H$ has a
$J$-skew-self-adjoint (respectively $J$-unitary) extension in a Hilbert space $\widetilde H\supseteq H$.
We shall follow the ideas of Galindo in~\cite{cit_450_G} with necessary modifications.
In particular, Lemma in~\cite{cit_450_G} can not be applied in our case, since its assumptions can never be satisfied with
$T$: $T^2 = I$. In fact, in this case $T$ would be a conjugation in $H$. Choosing an element $f\in H$ of an orthonormal basis
in $H$ which corresponds to $T$,
we would get $(f,Tf)= (f,f) = 1\not= 0$. Moreover, an exit out of the original space can appear in our case.

\noindent
We notice that under stronger assumptions on a $J$-skew-symmetric operator the existence of
a $J$-skew-self-adjoint extension was proved by Kalinina in~\cite{cit_550_K}.

\noindent
\textbf{Notations.}
As usual, we denote by $\mathbb{R}, \mathbb{C}, \mathbb{N}, \mathbb{Z}, \mathbb{Z}_+$,
the sets of real numbers, complex numbers, positive integers, integers and non-negative integers,
respectively. Set $\overline{0,d} = \{ 0,1,...,d \}$, if $d\in \mathbb{N}$; $\overline{0,\infty} = \mathbb{Z}_+$.
If H is a Hilbert space then $(\cdot,\cdot)_H$ and $\| \cdot \|_H$ mean
the scalar product and the norm in $H$, respectively.
Indices may be omitted in obvious cases.
For a linear operator $A$ in $H$, we denote by $D(A)$
its  domain, by $R(A)$ its range, and $A^*$ means the adjoint operator
if it exists. If $A$ is invertible then $A^{-1}$ means its
inverse.
For a set $M\subseteq H$
we denote by $\overline{M}$ the closure of $M$ in the norm of $H$.
By $\mathop{\rm Lin}\nolimits M$ we denote
the set of all linear combinations of elements of $M$,
and $\mathop{\rm span}\nolimits M:= \overline{ \mathop{\rm Lin}\nolimits M }$.
By $E_H$ we denote the identity operator in $H$, i.e. $E_H x = x$,
$x\in H$. In obvious cases we may omit the index $H$.
All appearing Hilbert spaces are assumed to be separable.

\section{Extensions of $J$-skew-symmetric and $J$-isometric operators.}

We shall make use of the following lemma.

\begin{lem}
\label{l2_1}
Let $H$ be a Hilbert space with a positive even or infinite dimension, and $J$ be a conjugation on $H$.
Then there exists a subspace $M$ in $H$ such that
$$ M \oplus JM = H. $$
\end{lem}
\textbf{Proof. }
Let $\{ f_n \}_{n=0}^{2d+1}$ be an orthonormal basis in $H$ corresponding to $J$, i.e. such that $Jf_n = f_n$, $0\leq n\leq d$;
$d\in Z_+\cup\{ +\infty \}$ ($2d+2=\dim H$).
Set
$$ f_{2k,2k+1}^+ = \frac{1}{\sqrt{2}}( f_{2k} + if_{2k+1} ),\quad f_{2k,2k+1}^- = \frac{1}{\sqrt{2}}( f_{2k} - if_{2k+1} ),\qquad
k\in \overline{0,d}. $$
It is easy to see that $\{ f_{2k,2k+1}^+, f_{2k,2k+1}^- \}_{k=0}^d$ is an orthonormal basis in $H$.
Set $M := \mathop{\rm span}\nolimits \{ f_{2k,2k+1}^+ \}_{k=0}^d$. It remains to notice that
$JM = \mathop{\rm span}\nolimits \{ f_{2k,2k+1}^- \}_{k=0}^d$.
$\Box$

\begin{thm}
\label{t2_1}
Let $H$ be a Hilbert space and $J$ be a conjugation on $H$. Let $A$ be a $J$-skew-symmetric ($J$-isometric) operator
in $H$. Suppose that $\overline{D(A)}=H$ (respectively $\overline{D(A)}=\overline{R(A)}=H$).
Then there exists a $J$-skew-self-adjoint (respectively $J$-unitary) extension of $A$ in a Hilbert space $\widetilde H\supseteq H$
(with an extension of $J$ to a conjugation on $\widetilde H$).
\end{thm}
\textbf{Proof.} Let $A$ be such an operator as that in the statement of the theorem.
The operator $A$ admits the closure which is $J$-skew-symmetric (respectively $J$-isometric) (see, e.g.~\cite[p. 18]{cit_300_Z}).
Thus, without loss of generality we shall assume that $A$ is closed.
In what follows, in the case of a $J$-skew-symmetric ($J$-isometric) $A$, we shall say about
case~(a) (respectively case~(b)).
Set $H_2 = H\oplus H$, and consider the following  transformations on $H_2$:
$$ J_2 \{ x, y \} = \{ Jx, Jy \},\ V \{ x, y \} = \{ y, -x \},\ U \{ x, y \} = \{ y, x \},\quad \forall \{ x,y \}\in H_2, $$
and $R := UJ_2 = J_2U$, $K := VR$. Observe that $R$ and $K$ are conjugations on $H_2$.
The graph of an arbitrary linear operator $C$ in the Hilbert space $H$ will be denoted by $G_C$ ($\subseteq H_2$).
Observe that
\begin{equation}
\label{f2_5}
J_2 G_C = G_{JCJ},\quad RG_C = UG_{JCJ}.
\end{equation}
If $\overline{D(C)}=H$, then
\begin{equation}
\label{f2_6}
G_{C^*} = H_2\ominus VG_C.
\end{equation}
In the case~(a) we may write:
$$ (\{ x,Ax \}, \{ JAJ y,y \}) = (x,JAJy) + (Ax,y) = 0,\quad \forall x\in D(A), y\in D(JAJ). $$
Then
\begin{equation}
\label{f2_7}
G_A \perp RG_A.
\end{equation}
In the case~(b), we have
$$ (\{ x,Ax \}, \{ JA^{-1}J y,-y \}) =  0,\quad \forall x\in D(A), y\in D(JA^{-1}J), $$
and therefore
\begin{equation}
\label{f2_9}
G_A \perp KG_A.
\end{equation}
Set $D = \left\{ \begin{array}{cc} H_2\ominus [G_A\oplus RG_A] & \mbox{in the case (a)} \\
H_2\ominus [G_A\oplus KG_A] & \mbox{in the case (b)}\end{array}\right.$.
If $D=\{0\}$ then it means that $A$ is $J$-skew-self-adjoint (respectively $J$-unitary), see considerations for the operator $B$ below.
In the opposite case, we have $RD=D$ (respectively $KD=D$).

At first, suppose that $D$ has a positive even or infinite dimension.
By~Lemma~\ref{l2_1} we obtain that there exists a subspace $X\subseteq D$ such that
$X\oplus RX = D$ (respectively $X\oplus KX = D$).
Since each element of $X$ is orthogonal to $RG_A = VG_{-JAJ}$ ($KG_A = VG_{JA^{-1}J}$), by~(\ref{f2_6}) it follows that
\begin{equation}
\label{f2_11}
X\subseteq G_{-JA^*J} \quad (\mbox{respectively }X\subseteq G_{J(A^{-1})^*J}).
\end{equation}
Set $G' = G_A\oplus X$. Suppose that $\{ 0,y \}\in G'$. Then there exist $\{ x,Ax \}\in G_A$ such that $\{0,y\} -
\{ x,Ax \} = \{ -x,y-Ax\}\in X$. By~(\ref{f2_11}) we get
$y-Ax = JA^*Jx$ (respectively $y-Ax = -J(A^{-1})^*Jx$), and therefore $y=0$.
Thus, $G'$ is a graph $G_B$ of a densely defined linear operator $B$. Moreover, we have
$$ G_B\oplus RG_B = H_2 \quad (\mbox{respectively } G_B\oplus KG_B = H_2). $$
In the case~(a) we get
$$ UG_B\oplus URG_B = H_2; $$
$$ G_{ (-B)^*} = H_2 \ominus VG_{-B} = H_2 \ominus UG_B = URG_B = J_2 G_B = G_{JBJ}. $$

\noindent
In the case~(b) we get
$$ VG_B\oplus VKG_B = H_2; $$
$$ G_{B^*} = H_2 \ominus VG_B =  VKG_B = - RG_B = G_{JB^{-1}J}. $$

Suppose now that $D$ has a positive odd dimension.
In this case we consider a linear operator $\mathcal{A} = A\oplus A$, with $D(\mathcal{A}) = D(A)\oplus D(A)$, in a Hilbert space
$\mathcal{H} = H\oplus H$ with a conjugation $\mathcal{J} = J\oplus J$.
Observe that $\mathcal{A}$ is a closed $\mathcal{J}$-skew-symmetric ($\mathcal{J}$-isometric) operator with
$\overline{D(\mathcal{A})} = \mathcal{H}$ (respectively $\overline{D(\mathcal{A})} = \overline{R(\mathcal{A})} = \mathcal{H}$).
Its graph $G_{\mathcal{A}}$ in a Hilbert space $\mathcal{H}_2 = \mathcal{H}\oplus \mathcal{H}$ may be identified
with $G_A\oplus G_A$ in $H_2\oplus H_2$:
$$ G_{\mathcal{A}} = \{ \{ (f,Af), (g,Ag) \},\ f,g\in D(A) \}. $$
Let $\mathcal{R}$, $\mathcal{K}$ be constructed for $\mathcal{A}$ as $R$ and $K$ for $A$.
%Then $\mathcal{R} = R\oplus R$, $\mathcal{K} = K\oplus K$.
In the case (a) we see that
$$ \mathcal{H}_2\ominus [G_\mathcal{A}\oplus \mathcal{R} G_\mathcal{A}] =
(H_2\ominus [G_A\oplus RG_A]) \oplus (H_2\ominus [G_A\oplus RG_A]), $$
has a positive even dimension.
In the case (b), $\mathcal{H}_2\ominus [G_\mathcal{A}\oplus \mathcal{K} G_\mathcal{A}]$
has a positive even dimension.
Thus, we may apply the above construction with $\mathcal{A}$ instead of $A$.

\noindent
$\Box$

\begin{center}
{\large\bf On extensions of $J$-skew-symmetric and $J$-isometric operators.}
\end{center}
\begin{center}
{\bf S.M. Zagorodnyuk}
\end{center}

In this paper it is
proved that each densely defined $J$-skew-symmetric operator (or each $J$-isometric operator with $\overline{D(A)}=\overline{R(A)}=H$) in a Hilbert space $H$
has a $J$-skew-self-adjoint (respectively $J$-unitary) extension in a Hilbert space $\widetilde H\supseteq H$.
We follow the ideas of Galindo in~[A.~Galindo, On the existence of $J$-self-adjoint extensions of $J$-symmetric operators with adjoint,
Communications on pure and applied mathematics, Vol. XV, 423-425 (1962)] with necessary modifications.
}

\end{document}